\documentclass[12pt]{article}
\usepackage[numbers,sort&compress]{natbib}
\usepackage{amssymb,amsmath,mathrsfs,amsfonts,amsthm}
\usepackage{enumerate}
\usepackage{paralist}
\usepackage{color}
\usepackage{epsfig}
\usepackage{graphicx}
\usepackage{setspace}
\usepackage{appendix}
\begin{document}
 \def\pd#1#2{\frac{\partial#1}{\partial#2}}
\def\dfrac{\displaystyle\frac}
\let\oldsection\section
\renewcommand\section{\setcounter{equation}{0}\oldsection}
\renewcommand\thesection{\arabic{section}}
\renewcommand\theequation{\thesection.\arabic{equation}}

\newtheorem{thm}{Theorem}[section]
\newtheorem{cor}[thm]{Corollary}
\newtheorem{lem}[thm]{Lemma}
\newtheorem{prop}[thm]{Proposition}
\newtheorem*{con}{Conjucture}
\newtheorem*{questionA}{Question}
\newtheorem*{thmA}{Theorem A}
\newtheorem*{thmB}{Theorem B}
\newtheorem{remark}{Remark}[section]
\newtheorem{definition}{Definition}[section]

\title{Global dynamics of  a parabolic type equation arising from the curvature flow
\thanks{The first author is supported by Chinese NSF (No. 11501207).
The second author is  supported by  NSF  of China (No. 11431005).
The third author is  supported by  NSF  of China (No. 11871148).} }

\author{Xueli Bai{\thanks{E-mail: mybxl110@163.com}}\\ {Department of Applied Mathematics, Northwestern Polytechnical
University,}\\{\small 127 Youyi Road(West), Beilin 710072,
Xi'an, P. R. China.}\\ Fang Li{\thanks{Corresponding author.
E-mail: lifang55@mail.sysu.edu.cn}}\\ {School of Mathematics, Sun Yat-sen
University,}\\{\small No. 135, Xingang Xi Road, Guangzhou 510275, P. R. China.} \\ Xiaoliu Wang {\thanks{
E-mail: xlwang@mail.seu.edu.cn}}\\ {School of Mathematics, Southeast University,}\\{\small No. 2, Sipailou, Nanjing 210096, P. R. China.} }

\date{}
\maketitle{}

\begin{abstract}
This paper studies a type of degenerate parabolic problem with nonlocal term
\begin{equation*}
  \begin{cases}
    u_t=u^p(u_{xx}+u-\bar{u})  & 0<t<T_{{\max}},\ 0<x<a,\\
    u_x(0,t)=u_x(a,t)=0     & 0<t<T_{{\max}},\\
    u(x,0)=u_0(x)   & 0<x<a,
  \end{cases}
\end{equation*}
where $p>1$, $a>0$. In this paper, the classification of the finite-time blowup/global existence phenomena based on the associated energy functional and explicit expression of all nonnegative steady states are demonstrated. Moreover, we combine the applications of Lojasiewicz-Simon inequality and energy estimates to derive that any bounded solution with positive initial data converges to some steady state as $t\rightarrow +\infty$.
\end{abstract}

{\bf Keywords}:
\vskip3mm {\bf MSC (2010)}: Primary:  ; Secondary:


\section{Introduction}
In this paper, we investigate global dynamics of solutions to the following parabolic type equation
\begin{equation}\label{mainpb}
  \begin{cases}
    u_t=u^p(u_{xx}+u-\bar{u})  & 0<t<T_{{\max}}, 0<x<a,\\
    u_x(t,0)=u_x(t,a)=0     & 0<t<T_{{\max}},\\
    u(0,x)=u_0(x)   & 0<x<a,
  \end{cases}
\end{equation}
where $p>1$, $\bar{u}=\frac{1}{a}\int_0^a u\,{\rm d}x$ and $T_{\max}\leq +\infty$ denotes the maximal time.
Throughout this paper, we always assume that the initial data $u_0$ satisfies
$$
u_0\in C([0,a]) \  \textrm{and}\ u_0>0 \ \ \textrm{in}\ [0,a].
$$

A typical property of the problem (\ref{mainpb}) is that the nonlocal term $\bar{u}$ in the equation keeps the following integral conserved:
  \begin{eqnarray}
  \int_0^a u^{1-p}(t,x)\,{\rm d}x=\int_0^a u_0^{1-p}(x)\,{\rm d}x   \  \ \
  \forall \  \ 0\leq t<T_{{\max}}.
  \label{integral}\end{eqnarray}
Moreover, due to the special structure of nonlocal term, the comparison principle does not always hold for the problem  (\ref{mainpb}). However, thanks to \cite[Lemma 2.1]{LWW}, the positivity of initial data guarantees the positivity of the corresponding solution for $0<t<T_{\max}$. Hence, throughout this paper, our discussion is restricted to positive solutions and nonnegative steady states of the problem (\ref{mainpb}).

Indeed, the study motivation of (\ref{mainpb}) arises from the popular topic of curvature flow. If $a=m\pi (m\in{\mathbb{N}^+})$ and the initial data $u_0$ satisfy
$$\int_0^a\frac{\cos x}{u_0^{p-1}}\,{\rm d}x = 0,$$
  it can be easily seen that the even-extension of   the  function $\kappa:=u_0^{p-1}$ could make $\kappa$ become a periodic function, which is nothing but the curvature function of
 some  closed, locally convex curve $\gamma_0$. Note that $\gamma_0$ is axis-symmetric and has $2m\pi$ total curvature. Then the problem  (\ref{mainpb}) could be used to describe the motion of $\gamma_0$ driven by a length-preserving curvature flow, see related studies in \cite{TW1,WLC,STW} and references therein. In particular, when $a=\pi$, $\gamma_0$ is a simple closed convex curve and the flow converges to a round circle as time goes to infinity according to the result in \cite{TW1}.

 One can also regard the problem  (\ref{mainpb}) as a natural generalization of the following degenerate parabolic problem not in divergence form
\begin{equation}\label{degepb}
  \begin{cases}
    v_t=v^p(v_{xx}+v)   & 0<t<T_{\max}, \ 0<x<a,\\
    v(0,t)=v(a,t)=0   & 0<t<T_{\max},\\
    v(x,0)=v_0(x)    & 0<x<a.
  \end{cases}
\end{equation}
Here, $p>0$ and $v_0\in C([0,a])$ is positive in the open interval $(0,a)$.  The problem   (\ref{degepb}) arises in the model for the resistive diffusion of a force-free magnetic field in a plasma confined between two walls, see \cite{Low}. It is well-known that the size of $[0,a]$ plays an important role in global solvability of solution $v$, see  \cite{FM,WWX,Win1,Win2,Win3} and references therein. Roughly speaking, related results can be summarized as follows:
\begin{itemize}
\item If $a>\pi$, then every solution $v$ blows up in  finite time.
 In addition, when $0<p<2$, its blowup rate is exactly $(T_{\max}-t)^{-1/p}$;
 when $p\geq 2$, blowup occurs essentially faster than the rate
 $(T_{\max}-t)^{-1/p}$.
 \item  If $a=\pi$, then all solutions are global in time.
 Moreover, when $0<p<3$, all solutions are bounded and converge to steady states as $t\to+\infty$;
 when $p\geq 3$, there are both bounded solution and unbounded solution,
 which hinges on $v_0(x)$.
 \item If $0<a<\pi$, then each solution exists globally
 and tends to zero as $t\to+\infty$.
\end{itemize}

Compared with the above problem,  the problem  (\ref{mainpb}) has an additional spatially nonlocal term $\bar{u}$. We are interested in illustrating how this modification influences the dynamics of solutions. An exploring study in \cite{LWW} shows that unlike the conclusion for the problem (\ref{degepb}),  the behavior of solutions to the problem (\ref{mainpb}) is not only influenced by the size of domain, but also impacted by positivity of the  associated energy  functional defined as
 $$
E[u](t) = \int_0^a \left[u_x^2 - (u-\bar u)^2  \right]\, {\rm d}x.
$$
To be more specific, the following results are established in  \cite{LWW}.
\begin{thmA}[\cite{LWW}]
Let $E[u](0)<0$. Then the solution of the  problem (\ref{mainpb}) blows up at  the maximal existence time $T_{\max}\leq +\infty$, and in particular, $T_{\max}<+\infty$ when $p>2$.
\end{thmA}

\begin{thmB}[\cite{LWW}]
Assume $a\leq \pi$ and $1<p<2$, or $a<\pi$ and  $p\geq 2$. Then, every solution $u$  of problem (\ref{mainpb}) is global and uniformly bounded. Moreover, $u$ converges smoothly to a constant for $a<\pi$ and $p>1$, and while, for $a=\pi$ and $1<p<2$, for any sequence $\{ t_j\}_{j=1}^{+\infty} \rightarrow +\infty$, there is a   subsequence   $\{ t_{jk}\}_{k=1}^{+\infty}$, such that $u(x, t_{jk})$ converges to the function $A\cos x +B$ in $C([0,a])$ as $k\rightarrow +\infty$, where $A$ and $B$ are constants satisfying $|A|\leq B$.
\end{thmB}

In this paper, we continue the studies in \cite{LWW} and provide a quite complete classification of the behavior of solutions to the problem (\ref{mainpb}).

First of all,  we manage to classify the  finite-time  blowup/global existence  behavior of solutions to the problem (\ref{mainpb}) based on the associated energy  functional $E[u](t)$.

\begin{thm}\label{thm-negativeenergy}
For $p>1$,  if $E[u](0)<0$, then the  solution $u(t,x)$  of the problem (\ref{mainpb})  blows up in finite time, i.e. $T_{\max}<+\infty$.
\end{thm}

This improves the result in Theorem A when $1<p\leq 2$.  Obviously, Theorem \ref{thm-negativeenergy} indicates that for any solution of the problem (\ref{mainpb}), if the associated energy  is negative sometime, then the solution blows up in finite time. Thus a  natural  question  is what happens if the associated energy $E[u](t)$ is always nonnegative.  Does the solution exist globally? If it does, is the solution bounded?  An affirmative answer is provided when $1<p \leq 2$, while it remains unknown when $p>2$.

\begin{thm}\label{thm-nonnegativeenergy}
Assume that $1<p\leq 2$, then one of the following alternatives must happen for the solution $u(t,x)$ of the problem (\ref{mainpb}):
\begin{itemize}
\item[(i)] if there exists $t_0\in [0,T_{\max})$ such that $E[u](t_0)<0$, then $T_{\max}<+\infty$;
\item[(ii)] if $E[u](t)\geq 0$ for $0\leq t<T_{\max}$, then  $T_{\max}=+\infty$ and  $u$ is bounded.
\end{itemize}
\end{thm}

Note that when $0<a \leq \pi$, according to Poincar\'{e} inequality, $E[\phi]\geq 0$ for any $0\not\equiv \phi\in H^1((0,a))$. Hence thanks to Theorem  \ref{thm-nonnegativeenergy}, if $1<p\leq 2$ and $0<a \leq \pi$, every solution of the problem (\ref{mainpb}) exists globally and is  bounded.  This improves the result in Theorem B, where the case $p=2, a=\pi$  is not covered.

Our next result is about  global convergence of bounded solutions to  the problem (\ref{mainpb}).

\begin{thm}\label{thm-convergence}
Assume that $u(t,x)$ is a bounded solution of the problem (\ref{mainpb}).  Then  $u(t,\cdot)$  converges  to a steady state $u^*(\cdot)$  in $C([0,a])$ as $t\rightarrow +\infty$.
\end{thm}

For the classical semilinear parabolic problems, it is known that all bounded solutions converge to some steady state in the one-dimensional case \cite{Matano, Z68}. However, when the nonlocal term is included, the dynamical phenomena become much more complicated \cite{FPolacik, Polacik}. Thus it is natural to investigate the sufficient conditions for the global convergence of bounded solutions. See for example \cite{HMatanoNW}.

However, in the problem (\ref{mainpb}),  the appearance of nonlocal term  results in the lack of comparison principles, while the nonlinearity $u^p$ in front of $u_{xx}$ makes the equation degenerate when the solution touches zero. Therefore, the problem (\ref{mainpb}) is dramatically different from the classical semilinear parabolic problems and it becomes  much more challenging  to demonstrate the global convergence of bounded solutions.

For simplicity, denote the $\omega-$limit set of the problem (\ref{mainpb}) as follows.
$$
\omega[u_0]= \omega[u_0 \, |\, C([0,a])] =  \bigcap_{\tau\geq 0} \overline{ \{ u(t,\cdot; u_0)\, |\, t\geq \tau \}},
$$
where assume that $u(t,x;u_0)$ is global and bounded.

To prove Theorem \ref{thm-convergence}, first observe that since the solution $u(t,x)$ is bounded,  $\omega[u_0]$ is not empty. Then it suffices to show that
$\omega[u_0]$  consists of a single element only. For this purpose,   different strategies are employed ingeniously based on the properties of $u^*$. To be more specific, fix any $u^*\in \omega[u_0]$.
\begin{itemize}
\item When $u^*>0$ in $[0,a]$, we introduce  a type of Lojasiewicz-Simon inequality designed for the problem (\ref{mainpb}) to prove the global convergence. 
\item When  $u^*$   touches  zeros somewhere,  we make use of the associated energy  functional  and some detailed analysis to verify the global convergence.  
\end{itemize}

It is worth pointing out that Theorem \ref{thm-convergence} automatically  indicates that for the problem (\ref{mainpb}),  $\omega-$limit set consists of a single element, which is a  steady state.  This is an important property itself.

Furthermore, thanks to Theorem \ref{thm-convergence}, we obtain the following byproduct.

\begin{prop}\label{prop}
Assume that $1<p\leq 2$ and $0<a \leq \pi$, or $p>2$ and $0<a<\pi$, then every solution of the problem (\ref{mainpb}) converges to some  nonnegative steady state  in $C([0,a])$ as $t\rightarrow +\infty$.
\end{prop}

This improves the corresponding result in Theorem B, where only subsequential convergence is obtained when  $1<p<2$ ,  $a =\pi$.

At the end, a complete classification and explicit expression of all nonnegative steady states to the problem (\ref{mainpb}) are  demonstrated in Propositions \ref{prop-steadystates} and \ref{prop-steadystates-initial}. In particular, Proposition \ref{prop-steadystates-initial}, together with Theorem \ref{thm-convergence}, reveals the role played by the length of the interval $a$ and the parameter $p$ in determining the global behavior of solutions to the problem (\ref{mainpb}).

This paper is organized as follows. Section 2 provides the proofs of Theorems \ref{thm-negativeenergy} and \ref{thm-nonnegativeenergy}.
Section 3 is devoted to  the  proofs of Theorem  \ref{thm-convergence} and Proposition \ref{prop}. In Section 4,  we present the statements and the proofs of Propositions \ref{prop-steadystates} and \ref{prop-steadystates-initial}. Some miscellaneous remarks are included in Section 5.

\section{Finite-time blowup/global existence}

This section is devoted to the proofs of Theorems \ref{thm-negativeenergy} and  \ref{thm-nonnegativeenergy}, where we try to investigate the relations between  the finite-time blowup/global existence phenomena in the  problem (\ref{mainpb})  and the  associated energy functional $E[u](t)$.

\begin{proof}[Proof of Theorem \ref{thm-negativeenergy}]
Suppose that  there exists $u_0(x)\in C([0,a])$ with  $u_0>0$ in $[0,a]$ and $E[u_0]<0$ such that the corresponding solution $u(t,x)$ exists globally. Note that
$$
{d\over \,{\rm d}t} E[u](t) = - 2\int_0^a u^{-p} u_t^2 \,{\rm d}x\leq 0
$$
and thus $E[u](t) \leq E[u_0]<0$.
\begin{itemize}
\item If $1<p<2$, let $y(t) = \int_0^a u^{2-p}(t,x) \,{\rm d}x$. Then
\begin{eqnarray*}
y'(t)  &=&  (2-p) \int_0^a u^{1-p} u_t \,{\rm d}x = -(2-p) E[u](t)\geq -(2-p) E[u_0]>0,\\
y''(t) &=& 2(2-p) \int_0^a u^{-p} u_t^2 \,{\rm d}x,
\end{eqnarray*}
and thus $\lim_{t\rightarrow +\infty} y(t)=+\infty$ and
$$
y''(t)y(t) =2(2-p) \int_0^a u^{-p} u_t^2 \,{\rm d}x \int_0^a u^{2-p}(t,x) \,{\rm d}x \geq {2\over 2-p} \left( y'(t)\right)^2.
$$
Then it is routine to see that
\begin{eqnarray*}
\left( y^{-{p\over 2-p}}\right)' &=& -{p\over 2-p} y^{-{2\over 2-p}} y',\\
\left( y^{-{p\over 2-p}}\right)''  &=& {p\over 2-p} y^{-{4-p \over 2-p}}\left[ {2\over 2-p} (y')^2 -y y'' \right]\leq 0.
\end{eqnarray*}
This is impossible since
$$
\lim_{t\rightarrow +\infty} y^{-{p\over 2-p}} =0.
$$
\item If $p = 2$, then
\begin{eqnarray*}
\left( -E[u](t) \right)^2  &=& \left(\int_0^a \left[u_x^2 - (u-\bar u)^2  \right] \,{\rm d}x\right)^2 =\left( \int_0^a  u^{-1}u_t \,{\rm d}x \right)^2\\
&\leq & a  \int_0^a u^{-2}u_t^2 \,{\rm d}x = {1\over 2}{d\over \,{\rm d}t} \left(-E[u](t)\right).
\end{eqnarray*}
Hence  $-E[u](t)$ blows up in finite time since $-E[u_0]>0$.
\item Assume that $p>2$.  Note that  $\int_0^a u^{1-p}(t,x) \,{\rm d}x = \int_0^a u_0^{1-p}(x) \,{\rm d}x$ since $\int_0^a u^{-p}u_t \,{\rm d}x = 0$. Then
\begin{eqnarray*}
\left( -E[u](t) \right)^2  &=& \left(\int_0^a \left[u_x^2 - (u-\bar u)^2  \right] \,{\rm d}x\right)^2 =\left( \int_0^a  u^{1-p}u_t \,{\rm d}x \right)^2\\
&\leq & \int_0^a u^{2-p} \,{\rm d}x  \int_0^a u^{-p}u_t^2 \,{\rm d}x = {1\over 2}\int_0^a u^{2-p} \,{\rm d}x {d\over \,{\rm d}t} \left(-E[u](t)\right)\\
&\leq & {1\over 2}a ^{1\over p-1}\left(\int_0^a u^{1-p} \,{\rm d}x \right)^{p-2\over p-1}{d\over \,{\rm d}t} \left(-E[u](t)\right)\\
&=& {1\over 2}a ^{1\over p-1}\left(\int_0^a u_0^{1-p} \,{\rm d}x \right)^{p-2\over p-1}{d\over \,{\rm d}t} \left(-E[u](t)\right)
\end{eqnarray*}
Thus, again  $-E[u](t)$ blows up in finite time since $-E[u_0]>0$.
\end{itemize}
The proof is complete.
\end{proof}

For the proof of  Theorem \ref{thm-nonnegativeenergy}, we prepare a useful lemma.

\begin{lem}\label{lem-bdd}
Let $u(t,x)$ denote a solution of the problem (\ref{mainpb}). If there exists $f(r)\in C([0,+\infty ))$ such that
$f(r)\geq 0$ in $[0,+\infty)$, $\lim_{r \rightarrow +\infty} f(r) =+\infty$ and  $\displaystyle\sup_{0<t<T} \int_0^a f(u(t,x))\,{\rm d}x< +\infty$, then $T=+\infty$ and $u$ is  bounded.
\end{lem}

\begin{proof}
Denote
$$
C_1 =   \sup_{0\not\equiv v \in H^1((0,a))} \frac{\| v\|_{L^{\infty}((0,a))}}{\| v\|_{H^1((0,a))}}  , \  C_2 = \sup_{0<t<T} \int_0^a f(u(t,x))\,{\rm d}x.
$$
Then one has
\begin{eqnarray*}
\| u\|^2_{L^{\infty}((0,a))} &\leq & C_1^2 \int_0^a \left(u^2 +u_x^2 \right)\,{\rm d}x\\
 &=& C_1^2  \left[\int_0^a \left(u^2 +(u-\bar u)^2 \right)\,{\rm d}x  + E[u](t)\right]\\
 &\leq & 2C_1^2 \int_0^a u^2 \,{\rm d}x + C_1^2 E[u_0]\\
 &\leq & 2C_1^2 \int_{\{u\leq R\}}u^2 \,{\rm d}x + 2C_1^2 \int_{\{u>R\}} {f(u)\over f(u)} u^2 \,{\rm d}x + C_1^2 E[u_0]\\
 &\leq & 2C_1^2 aR^2 + \frac{2C_1^2 C_2}{\inf_{r\in [R,+\infty)}f(r)}\| u\|^2_{L^{\infty}((0,a))}+ C_1^2 E[u_0].
\end{eqnarray*}
By choosing $R$ large such that
$$
\inf_{r\in [R,+\infty)}f(r) \geq 4C_1^2 C_2,
$$
it follows that
$$
\| u\|^2_{L^{\infty}((0,a))}  \leq  4C_1^2 aR^2 + 2 C_1^2 E[u_0].
$$
Therefore, $T=+\infty$ and  $u$ is globally bounded.
\end{proof}

\begin{proof}[Proof of Theorem \ref{thm-nonnegativeenergy}]
We discuss two cases $1<p<2$ and $p=2$ respectively.
\begin{itemize}
\item For $1<p<2$, choose $f(r) =  r^{2-p}$. Then
$$
{d\over \,{\rm d}t}\int_0^a f(u(t,x))\,{\rm d}x = \int_0^a (2-p)u^{1-p} u_t \,{\rm d}t=-(2-p) E[u](t)\leq 0.
$$
Hence by Lemma \ref{lem-bdd},  the desired conclusion follows for $1<p<2$.
\item For $p=2$, set $f(r) = \ln r$. Then
$$
{d\over \,{\rm d}t}\int_0^a f(u(t,x))\,{\rm d}x = \int_0^a u^{-1} u_t \,{\rm d}t=-E[u](t)\leq 0.
$$
Again Lemma \ref{lem-bdd} yields the desired conclusion for $p=2$.
\end{itemize}
The proof is complete.
\end{proof}

\section{Global convergence}
The purpose of this section is to demonstrate the proofs of Theorem \ref{thm-convergence} and its byproduct-Proposition \ref{prop}.

Let $u(t,x)$ denote the solution of the problem (\ref{mainpb}) and assume that $u(t,x)$  is   bounded for $t>0$, i.e., there exists $M>0$ such that $\| u(t, \cdot)\|_{C([0,a])}<M$  for $t>0$. Then
\begin{eqnarray}\label{H1-bdd}
\| u\|^2_{H^1((0,a))} &= & \int_0^a \left(u^2 +u_x^2 \right)\,{\rm d}x =  \int_0^a \left(u^2 +(u-\bar u)^2 \right)\,{\rm d}x  + E[u](t) \cr
 &\leq & 2\int_0^a u^2 \,{\rm d}x + E[u_0]< 2aM^2+ E[u_0].
\end{eqnarray}
Thus there exist a sequence $\{t_n \}_{n\geq 1}$ with $\lim_{n\rightarrow +\infty} t_n = +\infty$ and $u^* \in H^1((0,a))$ such that
\begin{equation}\label{pf-thm-H1}
u(t_n, \cdot) \rightharpoonup u^* \ \  \textrm{in} \ H^1((0,a))\ \ \textrm{as} \ n\rightarrow +\infty.
\end{equation}
Obviously, $u^*(x)\geq 0$ in $[0,a]$.

To prove Theorem \ref{thm-convergence},  it suffices to  show that
\begin{equation}\label{pf-thm-aim}
\lim_{t\rightarrow +\infty}  \|u(t,\cdot) - u^*(\cdot)\|_{C([0,a])} =0.
\end{equation}
For this purpose, two cases will be discussed  respectively since different approaches are employed.
\begin{itemize}
\item[{\it Case 1}.]  $u^*(x)>0$ in $[0,a]$.
\item[{\it Case 2}.] $u^* (x)$ admits at lease one zero in $[0,a]$.
\end{itemize}
Moreover, in {\it Case 2}, $u^*(x)$ will be expressed explicitly.
The proof of (\ref{pf-thm-aim}) is quite lengthy. We will prepare a series of lemmas  first, then handle  {\it Case 1},   {\it Case 2}  respectively.

\bigskip


First, we introduce  a type of Lojasiewicz-Simon inequality designed for our problem. Different from the original  Lojasiewicz-Simon inequality in \cite{Simon} (also see \cite{HJ} for a generalized version) which is valid for these functions close to some steady state, while the following one is valid for all $u\in H^2(0,a)$ with $u_x(0)=u_x(a)=0$. A short proof is included here for completeness.

\begin{lem}\label{lem-Simon}
For any $u\in H^2((0,a))$ with $u'(0)= u'(a)=0$, we have
\begin{equation}\label{Simon}
\int_0^a \left[u_x^2 -(u-\bar u)^2 \right]  \,{\rm d}x \leq C_a \int_0^a \left[  u_{xx} + u - \bar u \right]^2 \,{\rm d}x,
\end{equation}
where $C_a>0$ depends on $a$ only.
\end{lem}

\begin{proof}
Note that the eigenvalues of $-\Delta$ in $(0,a)$ with homogeneous Neumann boundary condition are $k^2\pi^2/a^2$,  with the corresponding normalized eigenfunctions
$$
\phi_k=\frac{\cos (k\pi x/a)}{\| \cos (k\pi x/a) \|_{L^2((0,a))}}, \ k\geq 0.
$$
Now assume that $u= \sum_{k=0}^{+\infty} c_k \phi_k = \bar u + \sum_{k=1}^{+\infty} c_k \phi_k$, where $c_k= \int_0^a u\phi_k \,{\rm d}x$, $k\geq 0$. Then it is standard to verify that
\begin{eqnarray*}
&& \int_0^a \left[u_x^2 -(u-\bar u)^2 \right]  \,{\rm d}x =  \sum_{k=1}^{+\infty} c_k^2 \left( {k^2\pi^2 \over a^2} -1 \right),\\
&& \int_0^a \left[  u_{xx} + u - \bar u \right]^2 \,{\rm d}x =  \sum_{k=1}^{+\infty}c_k^2 \left( {k^2\pi^2 \over a^2} -1 \right)^2.
\end{eqnarray*}
Thus obviously, there exists $C_a>0$, which depends on $a$ only, such that (\ref{Simon}) holds.
\end{proof}

Then we state the Aubin-Lions Lemma and Ascoli-Arzela Theorem which will be useful later.

\begin{lem}[Aubin-Lions Lemma]
Let $X_0,~X$ and $X_1$ be three Banach spaces and suppose that $X_0$ is compactly embedded in $X$ and that $X$ is continuously embedded in $X_1$, i.e., $X_0\hookrightarrow\hookrightarrow X\hookrightarrow X_1$. For $q >1$, let

$$W=\left\{u\in L^{\infty}((0,T);X_{0})\ |\ {\dot {u}}\in L^{q}((0,T);X_{1})\right\},$$
  then the embedding of $W$ into $C([0, T]; X)$ is compact.
\end{lem}

\begin{lem}[Ascoli-Arzela Theorem]
  Let $X$ be a compact Hausdorff space. Then a subset $F$ of $C(X)$ is relatively compact in the topology induced by the uniform norm if and only if it is equicontinuous and pointwise bounded.
\end{lem}

\bigskip

Now we consider {\it Case 1}, i.e., $u^*(x)>0$ in $[0,a]$.
Denote $\sigma= \min_{x\in [0,a]} u^*(x)>0$ and $u_n(t,x) = u(t+ t_n, x)$.
For clarity, some of key properties  will be stated as lemmas.

\begin{lem}\label{lem-precompact}
$\{ u_n\}_{n\geq 1}$ is precompact in $C([0,1]\times [0,a])$.
\end{lem}
\begin{proof}
 By (\ref{H1-bdd}), one has
    $$
    \sup_{n\geq 1} \| u_n\|_{L^{\infty}\left((0,1); H^1((0,a))\right)}< +\infty.
    $$
    Also
    \begin{eqnarray*}
    E[u](t_n)- E[u](t_n+1)& =&  2\int_{t_n}^{t_n+1}\int_0^a u^{-p} u_t^2 \,{\rm d}x {\rm d}t\geq 2M^{-p}\int_{t_n}^{t_n+1}\int_0^a  u_t^2 \,{\rm d}x {\rm d}t \\
    &\geq & 2M^{-p}\left\|  (u_n)_t \right\|^2_{L^{2}\left((0,1); H^1((0,a))\right)}.
    \end{eqnarray*}
    Thanks to Theorem \ref{thm-negativeenergy}, one sees that $E[u](t)\geq 0$ for $t>0$. Thus it follows that
    $$
    \left\|  (u_n)_t \right\|^2_{L^{2}\left((0,1); H^1((0,a))\right)} \leq {1\over 2} M^p E[u](t_n)\leq {1\over 2} M^p E[u_0].
    $$
    Then  by choosing $X_0 =H^1((0,a))$, $X= C([0,a])$ and $X_1 = L^2((0,a))$, Aubin-Lions Lemma can be applied to indicate that $\{ u_n\}_{n\geq 1}$ is precompact  $C([0, T]; X)$, i.e., $\{ u_n\}_{n\geq 1}$ is precompact in $C([0,1]\times [0,a])$.
\end{proof}

\begin{lem}\label{lem-lowerbound}
There exist a subsequence of $\{t_n\}_{n\geq 1}$, still denoted by $\{t_n\}_{n\geq 1}$, $\delta \in (0,1)$ and $N>0$ such that for $n\geq N$, $u(t_n+t,x)\geq {\sigma\over 2}$ in $[0,\delta] \times [0,a]$.
\end{lem}

\begin{proof}
  Since $ u_n(0, \cdot)= u(t_n, \cdot) \rightharpoonup u^*(\cdot) $ in $H^1((0,a))$ as $n\rightarrow +\infty$ and $H^1((0,a))\hookrightarrow\hookrightarrow  C([0,a])$, there exists  a subsequence of $\{t_n\}_{n\geq 1}$, still denoted by $\{t_n\}_{n\geq 1}$, such that
    $$
    u_n(0,\cdot) = u(t_n, \cdot) \rightarrow  u^*(\cdot) \  \textrm{in} \  C([0,a]) \ \ \textrm{as} \ n\rightarrow +\infty,
    $$
    and thus there exists $N>0$ such that for $n\geq N$, $u_n(0,x)> 3\sigma/4$ in $[0,a]$. Moreover, thanks to Lemma \ref{lem-precompact} and Ascoli-Arzela Theorem, the sequence $\{u_n\}_{n\geq 1}$ is equicontinuous in $[0,1]\times [0,a]$. Thus there exists $\delta\in (0,1)$ such that
    $$
    |u_n(t,x) - u_n(s,x)| <{\sigma\over 4}\ \ \textrm{for}\ |t-s|\leq \delta,  \ x \in [0,a],
    $$
    which implies that, for $n\geq N$,  $t\in [0,\delta]$, $x\in [0,a]$,
    \begin{eqnarray*}
    u(t_n+t, x)= u_n(t,x) \geq u_n(0,x) - |u_n(0,x) - u_n(t,x)|> {3\sigma\over 4}- {\sigma \over 4} = {\sigma\over 2}.
    \end{eqnarray*}
    The lemma is proved.
\end{proof}

Now  we are ready to  estimate $\| u(t,\cdot)\|_{H^2((0,a))}$ for certain time intervals.
According to the equation satisfied by $u$, one has for $n\geq N$
    \begin{eqnarray*}
    && \int_{t_n}^{t_n + {\delta\over 2}}\int_0^a u_{xx}^2 \,{\rm d}x {\rm d}t = \int_{t_n}^{t_n + {\delta\over 2}}\int_0^a \left[ u^{-p} u_t -u +\bar u \right]^2 \,{\rm d}x {\rm d}t\\
    &\leq & 2 \int_{t_n}^{t_n + {\delta\over 2}}\int_0^a \left( u^{-2p} u_t^2 +(u-\bar u)^2\right) \,{\rm d}x{\rm d}t  \leq 2 \left({\sigma\over 2}\right)^{-p} \int_{t_n}^{t_n + {\delta\over 2}}\int_0^a  u^{-p} u_t^2  \,{\rm d}x{\rm d}t + a \delta M^2\\
    &=& \left({\sigma\over 2}\right)^{-p} \int_{t_n}^{t_n + {\delta\over 2}} \left(- {d\over \,{\rm d}t} E[u](t)\right)  \,{\rm d}t + a \delta M^2\\
    &\leq & \left({\sigma\over 2}\right)^{-p}  E[u](t_n) + a \delta M^2\leq \left({\sigma\over 2}\right)^{-p}  E[u_0] + a \delta M^2,
    \end{eqnarray*}
     since $E[u](t)$ is decreasing  and $E[u](t)\geq 0$ for $t>0$ by Theorem \ref{thm-negativeenergy}. This indicates that there exists $s_n \in [t_n, t_n+ {\delta\over 2}]$ such that
     \begin{equation}\label{pf-thm-uxx-sn}
     \| u_{xx}(s_n, \cdot)\|_{L^2((0,a))} < +\infty\ \ \textrm{uniformly in}\ n.
     \end{equation}
     Hence by  Lemma \ref{lem-lowerbound}, (\ref{pf-thm-uxx-sn}) and parabolic regularity, it follows immediately that
     \begin{equation}\label{pf-thm-H2-n}
     \sup_{t\in [s_n, t_n +\delta]}\| u(t,\cdot)\|_{H^2((0,a))}\ \textrm{is uniformly bounded in}\ n.
     \end{equation}

The next result is crucial for the application of a type of Lojasiewicz-Simon inequality obtained in Lemma \ref{lem-Simon}.

\begin{lem}\label{lem-1-alpha}
There exist  $w \in C^2([0,a])$, $0<\alpha<{1\over 2}$ and a sequence $\{ r_n \}_{n\geq N}$ such that $w$ is a positive steady state of the problem (\ref{mainpb}), $w(x)\geq {\sigma \over 2}$ in $[0,a]$,  and $u(r_n,\cdot) \rightarrow w(\cdot)$ in $C^{1,\alpha}([0,a])$ as $n\rightarrow +\infty$.
\end{lem}

\begin{proof}
Notice that  for $n\geq N$
    \begin{eqnarray*}
    2\left({\sigma\over 2}\right)^{p} \int_{s_n}^{t_n + \delta }\int_0^a  u^{-2p} u_t^2  \,{\rm d}x {\rm d}t \leq 2 \int_{s_n}^{t_n + \delta }\int_0^a  u^{-p} u_t^2  \,{\rm d}x {\rm d}t = E[u](s_n) - E[u](t_n + \delta),
    \end{eqnarray*}
    which, together with the fact that $s_n \in [t_n, t_n+ {\delta\over 2}]$,  implies that there exists $r_n \in [s_n, t_n+\delta]$ such that as $n\rightarrow +\infty$,
    \begin{equation}\label{pf-thm-rn}
    \delta\left({\sigma\over 2}\right)^{p}  \int_0^a  u^{-2p}(r_n,x) u_t^2(r_n,x)  \,{\rm d}x  \leq E[u](s_n) - E[u](t_n + \delta) \rightarrow 0.
    \end{equation}

Moreover, thanks to  (\ref{pf-thm-H2-n}), one sees that $\| u(r_n,\cdot)\|_{H^2((0,a))}$ is uniformly bounded in $n\geq N$. Hence there exist $w \in C^{1,\alpha}([0,a])$, $0<\alpha<{1\over 2}$, and a subsequence of $\{ r_n \}_{n\geq N}$, still denoted by $\{ r_n \}_{n\geq N}$,  such that
\begin{equation}\label{pf-thm-C-alpha}
u(r_n,\cdot) \rightarrow w (\cdot)  \ \ \textrm{in} \ C^{1,\alpha}([0,a])  \ \ \textrm{as} \ n\rightarrow +\infty,
\end{equation}
and $w(x)\geq \sigma/2$ in $[0,a]$ due to  Lemma \ref{lem-lowerbound}.

Now for any given  $\phi\in H^1((0,a))$,
    $$
    \int_0^a u^{-p} u_t \,{\rm d}x = \int_0^a \left( -u_x\phi_x  + (u-\bar u)\phi \right)\,{\rm d}x,
    $$
    where by choosing $t=r_n$ and letting $n\rightarrow +\infty$, it follows from (\ref{pf-thm-rn}) and (\ref{pf-thm-C-alpha}) that
    $$
    \int_0^a \left( -w_x\phi_x  + (w -\overline {w})\phi \right)\,{\rm d}x=0,
    $$
    and thus  $w$ is a positive steady state of the problem (\ref{mainpb}).
\end{proof}

Also we claim that
\begin{equation}\label{pf-thm-E=0}
\lim_{t \rightarrow +\infty} E[u](t) =0.
\end{equation}
Suppose that this claim is not true. Recall that by Theorem \ref{thm-negativeenergy}, we have $E[u](t)\geq 0$ for $t>0$. Thus there exists $A>0$ such that $\lim_{t \rightarrow +\infty}  E[u](t) = A$.  Notice that on the side,
    if $p>1$, $p\neq 2$, then
    $$
    {d\over \,{\rm d}t} \int_0^a u^{2-p} \,{\rm d}x = - (2-p) E[u](t).
    $$
    On the other side, if $p=2$, then
    $$
    {d\over \,{\rm d}t} \int_0^a \ln u \,{\rm d}x = -  E[u](t).
    $$
    In either case, it contradicts to (\ref{pf-thm-H1}) with $u^*(x)>0$ in $[0,a]$.

Let $v(t,x) = u(t,x)- w(x)$, where $w(x)$ is  obtained in Lemma \ref{lem-1-alpha}. Then direct computation yields that
\begin{eqnarray}\label{pf-Eu=Ev}
  E[v](t)   =   \int_0^a \left[v_x^2 - (v-\bar v)^2  \right] \,{\rm d}x = E[u](t).
\end{eqnarray}
Thus thanks to Lemma \ref{lem-1-alpha} and (\ref{pf-thm-E=0}),  for any $\epsilon>0$, there  exists $N_1\geq N>0$, where $N$ is designated in Lemma \ref{lem-lowerbound}, such that
for   $t\geq r_{N_1}$,
\begin{eqnarray}\label{pf-thm-LS}
\| v(r_{N_1},\cdot)\|_{H^1((0,a))} <{\epsilon\over 4}, \ C_a^{1\over 2}\left( {4M\over \sigma}\right)^p \sqrt{E[v](t)}<{\epsilon \over 4}, \ \sqrt{E[v](t)} <  {\epsilon\over 2}.
\end{eqnarray}
Denote
$$
\tilde t = \sup \left\{ t \ \big|\  \| v(s,\cdot)\|_{H^1((0,a))} \leq 2\epsilon,\ \textrm{for} \ s\in [r_{N_1}, t]  \right\}.
$$
To show the desired conclusion that  $u(t,\cdot)$ converges to  $w(\cdot)$  in $C([0,a])$ as $t\rightarrow +\infty$, it suffices to show $\tilde t= +\infty.$ Suppose that $\tilde t< +\infty$.

First for $t\in [r_{N_1}, \tilde t]$,
$$
\| v(t, \cdot)\|_{L^{\infty}((0,a))} \leq C_1 \| v(t,\cdot) \|_{H^1((0,a))}\leq 2C_1\epsilon <{\sigma\over 4}
$$
by choosing $\epsilon$ smaller if necessary,
where
$$
C_1 =   \sup_{0\not\equiv v \in H^1((0,a))} \frac{\| v\|_{L^{\infty}((0,a))}}{\| v\|_{H^1((0,a))}} .
$$
This, together with Lemma \ref{lem-1-alpha},  indicates that for $t\in [r_{N_1}, \tilde t]$,
\begin{equation}\label{pf-thm-u-lowerbound}
u(t,x) = w(x) +v(t,x) \geq w(x)- \| v(t, \cdot)\|_{L^{\infty}((0,a))}>{\sigma\over 4}.
\end{equation}

Then  by Lemma \ref{lem-Simon}, (\ref{pf-Eu=Ev}), (\ref{pf-thm-u-lowerbound}), we have
\begin{eqnarray*}
&& - {d\over \,{\rm d}t}  \sqrt{E[v](t)} = - {d\over \,{\rm d}t}  \sqrt{E[u](t)} = -{1\over 2} \frac{1}{\sqrt{E[u](t)} }{d\over \,{\rm d}t}  E[u](t)\\
&=& \left(\int_0^a \left[u_x^2 - (u-\bar u)^2  \right] \,{\rm d}x  \right)^{-{1\over 2}}\int_0^a u^{-p} u_t^2 \,{\rm d}x\\
&\geq & M^{-p}\left(C_a \int_0^a \left[  u_{xx} + u - \bar u \right]^2 \,{\rm d}x\right)^{-{1\over 2}}\int_0^a  u_t^2 \,{\rm d}x\\
&\geq & M^{-p} C_a^{-{1\over 2}}\left( \int_0^a \left(  u^{-p} u_t\right)^2 \,{\rm d}x\right)^{-{1\over 2}}\int_0^a  u_t^2 \,{\rm d}x\\
&\geq &   C_a^{-{1\over 2}} \left({\sigma\over 4M}\right)^p \| u_t(t, \cdot)\|_{L^2((0,a))}=  C_a^{-{1\over 2}} \left({\sigma\over 4 M}\right)^p \| v_t(t, \cdot) \|_{L^2((0,a))}.
\end{eqnarray*}
It follows from (\ref{pf-thm-LS})  that
$$
\int_{r_{N_1}}^{\tilde t}\| v_t(t, \cdot) \|_{L^2((0,a))} \,{\rm d}t\leq C_a^{1\over 2} \left({4M\over \sigma}\right)^p \sqrt{E[v](r_{N_1})} <{\epsilon \over 4},
$$
which, together with (\ref{pf-thm-LS}),  implies that
$$
\| v(\tilde t, \cdot) \|_{L^2((0,a))} \leq \| v(r_{N_1}, \cdot) \|_{L^2((0,a))} + \int_{r_{N_1}}^{\tilde t}\| v_t(t, \cdot) \|_{L^2((0,a))} \,{\rm d}t <{\epsilon \over 2}.
$$
and thus again due to (\ref{pf-thm-LS}),  we have
\begin{eqnarray*}
&& \| v(\tilde t, \cdot) \|_{H^1 ((0,a))} \leq  \| v_x(\tilde t, \cdot) \|_{L^2((0,a))} +\epsilon\\
  &=&  \left( E[v](\tilde t) + \int_0^a(v(\tilde t, x) -\bar v (\tilde t))^2  \,{\rm d}x  \right)^{1\over 2} +\epsilon \leq \sqrt{E[v](\tilde t)} + \| v(\tilde t, \cdot) \|_{L^2((0,a))} +\epsilon\\
  &<& {\epsilon \over 2} + {\epsilon \over 2} +\epsilon = 2 \epsilon.
\end{eqnarray*}
This is a contradiction to the definition of $\tilde t$.
Therefore, indeed $u^*(x) \equiv w(x)$ is a positive steady state of the problem (\ref{mainpb}) and
$$
\lim_{t \rightarrow +\infty} u(t,\cdot) = u^*(\cdot)\ \ \textrm{in}\ C([0,a]).
$$


\bigskip

Next we focus on {\it Case 2}, i.e., $u^* (x)$ admits at lease one zero in $[0,a]$.
Denote
$$
\Omega_+ = \left\{ x\in [0,a]\ \big| \ u^*(x)>0 \right\}.
$$
First of all, we establish the following result.

\begin{lem}
$u^*$ satisfies
$$
(u^*)'' + u^* - \overline{u^*} =0\ \ \textrm{in}\ \Omega_+,
$$
$(u^*)' (0) =0$ if $0\in \Omega_+$ and  $(u^*)' (a) =0$ if $a\in \Omega_+$.
\end{lem}
\begin{proof}
For any $\sigma>0$, denote
$$
\Omega_{\sigma} = \left\{ x\in [0,a]\ \big| \ u^*(x)>\sigma \right\}.
$$
Same as the proof of Lemma \ref{lem-lowerbound}, there exist a subsequence of $\{t_n\}_{n\geq 1}$, still denoted by $\{t_n\}_{n\geq 1}$, $\delta \in (0,1)$ and $N>0$ such that for $n\geq N$, $u(t_n+t,x)\geq {\sigma\over 2}$ in $[0,\delta] \times \bar\Omega_{\sigma}$.

Choose $\rho(t)  \in C_c^{\infty}((0,\delta))$ with $\int_0^{\delta} \rho \,{\rm d}t =1$ and $\phi \in C_c^{\infty}(\Omega_{\sigma})$. Then
\begin{eqnarray*}
&& - \int_{t_n}^{t_n +\delta} \int_{\Omega_{\sigma}} {1\over 1-p} u^{1-p} \rho' (t-t_n) \phi(x) \,{\rm d}x {\rm d}t =   \int_{t_n}^{t_n +\delta} \int_{\Omega_{\sigma}} u^{-p}u_t \rho (t-t_n) \phi(x) \,{\rm d}x {\rm d}t\\
&=&\int_{t_n}^{t_n +\delta} \int_{\Omega_{\sigma}} (u_{xx}+u-\bar{u}) \rho (t-t_n) \phi(x) \,{\rm d}x {\rm d}t\\
&=&\int_{t_n}^{t_n +\delta} \int_{\Omega_{\sigma}} \left[u  \phi''(x) +  (u-\bar{u})  \phi(x) \right]\rho (t-t_n) \,{\rm d}x {\rm d}t,
\end{eqnarray*}
and by setting $u_n(t,x) = u(t+ t_n, x)$, it becomes
\begin{eqnarray}\label{pf-thm-weaksol}
&& \int_{0}^{\delta} \int_{\Omega_{\sigma}} \left[u_n  \phi''(x) +  (u_n-\bar{u}_n)  \phi(x) \right]\rho (t) \,{\rm d}x {\rm d}t = - \int_{0}^{\delta} \int_{\Omega_{\sigma}} {1\over 1-p} u_n^{1-p} \rho' (t ) \phi(x) \,{\rm d}x {\rm d}t \cr
 &=& - \int_{0}^{\delta} \int_{\Omega_{\sigma} } {1\over 1-p} u_n^{1-p}(0,x) \rho' (t ) \phi(x) \,{\rm d}x {\rm d}t\cr
  && - \int_{0}^{\delta} \int_{\Omega_{\sigma} } {1\over 1-p}\left[ u_n^{1-p}(t,x) -u_n^{1-p}(0,x)  \right]\rho' (t ) \phi(x) \,{\rm d}x {\rm d}t\cr
 &=& - \int_{0}^{\delta} \int_{\Omega_{\sigma} } {1\over 1-p}\left[ u_n^{1-p}(t,x) -u_n^{1-p}(0,x)  \right]\rho' (t ) \phi(x) \,{\rm d}x {\rm d}t.
\end{eqnarray}

Let us estimate $ u_n^{1-p}(t,x) -u_n^{1-p}(0,x)$. Indeed, we claim that
$$
\lim_{n\rightarrow +\infty} \int_{\Omega_{\sigma}}\left[ u_n^{1-p}(t,x) -u_n^{1-p }(0,x)  \right]^2 \,{\rm d}x =0.
$$
Direct computation yields that for $n\geq N$, $t\in [0,\delta]$
\begin{eqnarray*}
&& \int_{\Omega_{\sigma}}\left[ u_n^{1-p}(t,x) -u_n^{1-p }(0,x)  \right]^2 \,{\rm d}x = \int_{\Omega_{\sigma}}\left[(1-p)
\int_0^t u_n^{-p}(s,x) {\partial \over\partial s}u_n (s,x) \,{\rm d}s \right]^2 \,{\rm d}x\\
&\leq & (1-p)^2 \left({2\over \sigma}\right)^{p} \int_{\Omega_{\sigma}} t \int_{t_n}^{t_n+ t} u^{-p}(s,x)   u_s^2 (s,x)  \,{\rm d}s {\rm d}x\\
&\leq & \delta (1-p)^2 \left({2\over \sigma}\right)^{p}  \int_{t_n}^{t_n+ t} -{1\over 2} {d\over \,{\rm d}s} E[u](s) \,{\rm d}s \leq    {\delta\over 2}(1-p)^2 \left({2\over \sigma}\right)^{p}   E[u](t_n).
\end{eqnarray*}
Hence the claim follows from (\ref{pf-thm-E=0}). Then in (\ref{pf-thm-weaksol}), by letting $n\rightarrow +\infty$, we have
$$
\int_{0}^{\delta} \int_{\Omega_{\sigma}} \left[u^*  \phi''(x) +  (u^* -\overline{u^*})  \phi(x) \right]\rho (t) \,{\rm d}x {\rm d}t=0.
$$
Since $\rho(t)  \in C_c^{\infty}((0,\delta))$ is arbitrary,  this is equivalent to
\begin{equation}\label{pf-thm-weaksol-1}
\int_{\Omega_{\sigma}} \left[u^*  \phi''(x) +  (u^* -\overline{u^*})  \phi(x) \right]  \,{\rm d}x  =0.
\end{equation}

Now let $(a_1, a_2)$ denote a connected interval in $\Omega_+$ with $u^*(a_1) = u^*(a_2)=0.$ Then for any fixed $\phi \in C_c^{\infty}((a_1, a_2))$, similar to the derivation of (\ref{pf-thm-weaksol-1}), one has
\begin{equation*}
 \int_{a_1}^{a_2} \left[u^*  \phi''(x) +  (u^* -\overline{u^*})  \phi(x) \right]  \,{\rm d}x  =0,
\end{equation*}
which becomes
\begin{equation*}
 \int_{a_1}^{a_2} \left[-(u^*)'  \phi'(x) +  (u^* -\overline{u^*})  \phi(x) \right]  \,{\rm d}x  =0,
\end{equation*}
since $u^* \in H^1((0,a))$. Thus by elliptic regularity, we have $u^*$ is a classical solution of
$$
\begin{cases}
(u^*)'' + u^* -\overline{u^*} =0 & x\in (a_1, a_2),\\
u^*(a_1)=0,\  u^*(a_2)=0.
\end{cases}
$$

Next assume that $0\in \Omega_+$ and thus there exists $a_3\in(0,a]$ such that $u^*>0$ in $[0,a_3)$ and $u^*(a_3)=0$.
We will verify that $(u^*)'(0)=0$.

For this purpose, first fix any $\phi\in C^{\infty}((0,a_3))\bigcap C([0,a_3])$ with $\textrm{supp}\ \phi \subseteq [0,a_3)$.  Similar to the derivation of (\ref{pf-thm-weaksol-1}), one has
\begin{equation*}
u^*(0) \phi'(0)  + \int_0^{a_3} \left[u^*  \phi''(x) +  (u^* -\overline{u^*})  \phi(x) \right]  \,{\rm d}x  =0.
\end{equation*}
Since $u^* \in H^1((0,a))$, this is equivalent to
$$
\int_0^{a_3} \left[- (u^*)'  \phi'(x) +  (u^* -\overline{u^*})  \phi(x) \right]  \,{\rm d}x  =0.
$$
For $x\in [-a_3, a_3]$, define $v(x) = u^*(|x|)$. Then it is easy to see, for any $\phi \in C_c^{\infty}((-a_3, a_3))$, one has
$$
\int_{-a_3}^{a_3} \left[- v'  \phi'(x) +  (v -\overline{u^*})  \phi(x) \right]  \,{\rm d}x  =0.
$$
Hence $v\in H_0^1((-a_3, a_3))$ is a weak solution to the problem
$$
\begin{cases}
v'' + v -\overline{u^*} =0 & x\in (-a_3, a_3),\\
v(-a_3) =v (a_3)=0.
\end{cases}
$$
By elliptic regularity, it follows that $v\in C^2([-a_3, a_3])$. According to the definition of $v$, we have
$$
\begin{cases}
(u^*)'' + u^* -\overline{u^*} =0 & x\in (0, a_3),\\
(u^*)'(0)=0,\  u^*(a_3)=0.
\end{cases}
$$

Similarly, we can show that  $(u^*)' (a) =0$ if $a\in \Omega_+$.
\end{proof}

Now let us consider the set
$$
\Omega_0 = \left\{ x\in [0,a]\ \big| \ u^*(x)= 0 \right\}.
$$
Obviously, in {\it Case 2}, $\Omega_0 \neq \emptyset$.
We claim that {\it the set $\Omega_0$  consists of discrete points.}

Assume that there exist two distinct points $x_1, x_2 \in \Omega_0$. First of all, we need rule out the possibility that
$$
u^*(x) \equiv 0 \ \ \textrm{for}\ x\in [x_1, x_2].
$$
This is obvious, since  according to the problem (\ref{mainpb}), one sees that
 $$
\int_0^a u^{1-p}(t,x) \,{\rm d}x = \int_0^a u_0^{1-p} (x)\,{\rm d}x.
$$
Thus there exists $x_0\in (x_1, x_2)$ such that $u^*(x_0)>0$. Denote
$$
y_1= \inf \{  y\ | \ u^*(x)>0\ \textrm{for any} \ x\in [y, x_0]  \},\  y_2= \sup \{  y\ | \ u^*(x)>0\ \textrm{for any}\   x\in    [x_0,y]  \}.
$$
Thus $x_1 \leq y_1< y_2 \leq x_2$.
Then let $w = u^*/\overline{u^*}$ and it satisfies
$$
\begin{cases}
w '' + w - 1 =0 & x\in (y_1, y_2),\\
w(y_1)=0,\  w(y_2)=0.
\end{cases}
$$
Let $y_0 \in (y_1, y_2)$ denote the point where
$$
w(y_0) = \sup_{y\in (y_1, y_2)} w(y) = M>0.
$$
Then obviously
\begin{equation}\label{pf-thm--w}
\psi(x) = (M-1)\cos (x-y_0) +1
\end{equation}
satisfies
$$
\begin{cases}
\psi '' + \psi - 1 =0 & x\in (y_1, y_2),\\
\psi(y_0) = M=w(y_0), \ \psi'(y_0) =0= w'(y_0).
\end{cases}
$$
Hence it follows that  $w\equiv \psi $ in $[y_1,y_2]$ and there are three  observations.
\begin{itemize}
\item The interval $[y_1, y_2]$ is symmetric with respect to $x=y_0$, i.e., there exists $\ell>0$ such that $[y_1, y_2 ] = [y_0- \ell, y_0 +\ell]$.
\item $M\geq 2$, otherwise for $M\in (0,2)$, $w$ has no zeroes.
\item  $\pi/2 < \ell \leq  \pi$ since for $x \in [y_0, y_0+\pi /2]$, $w (x) >0$ and $w(y_0 +\pi) =-M+2 \leq 0.$
\end{itemize}

Till now, one sees that  the set $\Omega_0$  consists of discrete points, denoted by
$$
0\leq y_1 < y_2 < ,..., <y_n\leq a,\ n\geq 1,
$$
and the distance between two adjacent points  is bigger than $\pi$.
Here also it is routine to check that
\begin{equation}\label{pf-thm-integration}
\int_{y_1}^{y_2}w(x) \,{\rm d}x =2 (M-1)\sin \ell +y_2 -y_1\geq y_2 -y_1.
\end{equation}

Furthermore, if $0 \neq y_1$, then $w = u^*/\overline{u^*}$ satisfies
$$
\begin{cases}
w '' + w - 1 =0 & x\in (0, y_1),\\
w(y_1)=0,\  w'(0)=0.
\end{cases}
$$
Similar to previous discussions, we derive that
\begin{equation}\label{pf-thm-0-w}
w(x)= (w(0)-1)\cos x +1 \ \textrm{for}\ x\in [0, y_1],
\end{equation}
$w(0)\geq 2$, $\pi/2 < y_1 \leq \pi$ and
\begin{equation}\label{pf-thm-integration-half}
\int_{0}^{y_1}w(x) \,{\rm d}x = (w(0) - 1)\sin y_1+y_1  \geq y_1.
\end{equation}
If $y_n \neq a$, similarly we have
\begin{equation}\label{pf-thm-a-w}
w(x)= (w(a)-1)\cos (x-a) +1 \ \textrm{for}\ x\in [y_n, a],
\end{equation}
$w(a)\geq 2$, $\pi/2 < a- y_n \leq \pi$ and
\begin{equation}\label{pf-thm-integration-half-a}
\int_{y_n}^{a}w(x) \,{\rm d}x = ( w(a) -1)\sin (a- y_n) +y_1  \geq y_1.
\end{equation}

At the end, let us classify the case that $\Omega_0$  is not  empty.   Recall that $w = u^*/\overline{u^*}$, hence  thanks to (\ref{pf-thm-integration}), (\ref{pf-thm-integration-half})  and (\ref{pf-thm-integration-half-a})
\begin{eqnarray*}
1&=& {1\over a} \int_0^a w(x )\,{\rm d}x\\
&= &{1\over a} \sum_{i=1}^{n-1} \int_{y_i}^{y_{i+1}} w(x) \,{\rm d}x \geq {1\over a} \sum_{i=1}^{n-1} (y_{i+1}  - y_i)=1,
\end{eqnarray*}
which indicates that $y_{i+1} - y_i= 2 \pi$ for $1\leq i\leq n$, $y_1 =\pi$ if $y_1\neq 0$, and $y_n = a-\pi $ if $y_n\neq a$. Furthermore, it is easy to verify that $M=2$ by (\ref{pf-thm--w}), $w(0)=2$ by (\ref{pf-thm-0-w}) if $y_1\neq 0$ and   $w(a) =2$ by (\ref{pf-thm-a-w}) if $y_n \neq a$.

Therefore, in {\it Case 2},   i.e., $u^* (x)$ admits at lease one zero in $[0,a]$, $a= k \pi$, $k\geq 1$ and
$$
\textrm{either} \ u^*(x) = \overline{u^*}  (\cos x +1)  \ \textrm{or}\ u^*(x) = \overline{u^*}(\cos (x- \pi) +1).
$$
Moreover,  according to the problem (\ref{mainpb}),  we have
$$
 \int_0^a u_0^{1-p} (x)\,{\rm d}x  = \int_0^a u^{1-p}(t,x) \,{\rm d}x = \int_0^a (u^*)^{1-p}(t,x) \,{\rm d}x,
$$
which implies that $1<p<3/2$ and
\begin{equation}\label{pf-thm-twoss-express}
\textrm{either}\  u^* (x) = A (\cos x+1)\  \textrm{or}\   u^* (x) = A (\cos (x- \pi) +1)
\end{equation}
where
$$
A = \left[ \frac{\int_0^{k\pi} (\cos x +1)^{1-p}(t,x) \,{\rm d}x}{\int_0^{k\pi} u_0^{1-p} (x)\,{\rm d}x}\right]^{1\over p-1}.
$$

In other words, {\it Case 2} is possible only when $[0, a] = [0, k\pi]$, $k\geq 1$ and $1<p<3/2$. Also when {\it Case 2} happens,  the set of steady states of the problem (\ref{mainpb}) consists of two  distinct steady states as designated above in (\ref{pf-thm-twoss-express}).

At the end, we prove the global convergence of the solution to the problem (\ref{mainpb}) when {\it Case 2} happens. To show
$$
\lim_{t \rightarrow +\infty} u(t,\cdot) = u^*(\cdot)\ \ \textrm{in}\ C([0,a]),
$$
it suffices to show that if there exist two sequence $\{ t_n\}_{n\geq 1}$ and $\{\tau_n\}_{n\geq 1}$ such that
$$
\lim_{n \rightarrow +\infty} u(t_n,\cdot) = u_1(\cdot), \ \lim_{n \rightarrow +\infty} u(\tau_n,\cdot) = u_2(\cdot) \  \ \textrm{in}\ C([0,a]),
$$
where $u_1$ and $u_2$ are steady states of the problem (\ref{mainpb}), then $u_1 = u_2$. Suppose that this is not true, which means the set of steady states of the problem (\ref{mainpb}) consists of $u_1$ and $u_2$ as expressed in (\ref{pf-thm-twoss-express}).
Then there exists another sequence $\{ s_n\}_{n\geq 1}$ such that
\begin{equation}\label{pf-thm-twoss}
\| u(s_n,\cdot) -  u_i(\cdot) \|_{C([0,a])} \geq A,\ \ i=1,2.
\end{equation}
Thanks to (\ref{H1-bdd}), there exist $\tilde u\in H^1((0,k\pi))$ and  a subsequence of $\{ s_n\}_{n\geq 1}$ , still denoted by $\{ s_n\}_{n\geq 1}$, such that
$$
\lim_{t \rightarrow +\infty} u(s_n,\cdot) = \tilde u(\cdot)\ \ \textrm{in}\ C([0,a]).
$$
This, together with (\ref{pf-thm-twoss}), implies that
$$
\| \tilde u(\cdot) -  u_i(\cdot) \|_{C([0,a])} \geq {1\over 3}A ,\ \ i=1,2.
$$
Hence according to the previous arguments, $\tilde u$ should be the third steady state. This is a contradiction.  The proof of Theorem \ref{thm-convergence} is complete.

\bigskip

At the end, we include the proof of Proposition \ref{prop}.
\begin{proof}[Proof of Proposition \ref{prop}]
According to Theorem \ref{thm-convergence}, it suffices to show that when $1<p\leq 2$ and $0<a \leq \pi$, or $p>2$ and $0<a<\pi$, the solution of the problem (\ref{mainpb}) is always bounded.

First when $0<a \leq \pi$, according to Poincar\'{e} inequality, $E[\phi]\geq 0$ for any $0\not\equiv \phi\in H^1((0,a))$. Hence thanks to Theorem  \ref{thm-nonnegativeenergy}, if $1<p\leq 2$ and $0<a \leq \pi$, every solution of the problem (\ref{mainpb})  is  bounded.

However, when $p>2$ and $0<a<\pi$, thanks to Theorem A, the solution of the problem (\ref{mainpb}) is always bounded.
\end{proof}

\section{Nonnegative steady states}

First, we classify all the nonnegative steady states to the problem (\ref{mainpb}).
\begin{prop}\label{prop-steadystates}
Assume that $p>1$.
\begin{itemize}
\item[(i)] If $a\neq  k\pi$ for any $k\in\mathbb N$, then  the nonnegative steady states to the problem (\ref{mainpb}) consist of nonnegative constants only.
\item[(ii)]  If $a =  k\pi$ for some $k\in\mathbb N$, then  the nonnegative steady states to the problem (\ref{mainpb}) consist of $A\cos x+ B$ with $0\leq |A|\leq B$.
\end{itemize}
\end{prop}

\begin{proof}
First of all, by observation, one sees that  all nonnegative constants must be nonnegative steady states of (\ref{mainpb}).

Next, we discuss the existence/nonexistence of nonconstant and nonnegative steady states of the problem (\ref{mainpb}). Assume that $u^*\in C^2([0,a])$ is a nonconstant and nonnegative steady state to the problem (\ref{mainpb}). We claim that {\it $a =  k\pi$ for some $k\in\mathbb N$ and  $u^* = A\cos x+ B$ for some $A,B\in  \mathbb R$  with $0 < |A|\leq B$.}

Denote
$$
\Omega'_0 = \left\{ x\in [0,a]\ \big| \ (u^*)'(x)= 0 \right\}.
$$
It is easy to see that $\Omega'_0$ is discrete, denoted by
$$
0 =  y_1 < y_2 < ,..., <y_n = a,\ n\geq 1.
$$

Set $w = u^*/\overline{u^*}$ and it satisfies
$$
\begin{cases}
w '' + w - 1 =0 & x\in (0, a),\\
w'(0)=0,\  w'(a)=0.
\end{cases}
$$
Let $y_i \in [0, a]$, for some  $1\leq i\leq n-1$,  denote the point where
$$
w(y_i) = \sup_{x\in (0, a)} w (x) = M>0.
$$
Consider the function
\begin{equation*}
\psi(x) = (M-1)\cos (x-y_i) +1,
\end{equation*}
which satisfies
$$
\begin{cases}
\psi '' + \psi - 1 =0 & x\in (0, a),\\
\psi(y_i) = M=w(y_i), \ \psi'(y_i) =0= w'(y_i).
\end{cases}
$$
Hence it follows that  $w\equiv \psi $ in $[0,a]$. Moreover it is standard to deduce the following properties.
\begin{itemize}
\item $a =  k\pi$ for some $k\in\mathbb N$, $\Omega'_0 = \{ j\pi,\ 0\leq j\leq k \}$;
\item $M \in   (1,2]$ since the solution is nonnegative, nonconstant and $M=\sup_{x\in (0, a)} w (x) $;
\item  $w = (M-1)\cos x +1$ or $w = -(M-1)\cos x +1$.
\end{itemize}
These properties imply that
$$
u^* = A\cos x+ B \ \ \ \textrm{in}\ [0,k\pi]
$$
for some  $A,B\in  \mathbb R$   with  $0 < |A|\leq B$.

Therefore, nonconstant steady states of the problem (\ref{mainpb}) exist only when the size of the domain is a multiple of $\pi$.
\end{proof}

Furthermore, thanks to Theorem \ref{thm-convergence}, for any given positive initial data $u_0\in C([0,a])$,  the corresponding solution $u(t,x)$, if bounded,  converges to some nonnegative steady state, denoted by $u^*(x)$, i.e., $ \omega[u_0]= \{ u^* \}$. The possible explicit expression for $u^*$ is  provided as follows.

\begin{prop}\label{prop-steadystates-initial}
Assume that $p>1$, $u_0\in C([0,a])$ is a positive initial data to the problem (\ref{mainpb}) and $ \omega[u_0]= \{ u^* \}$.
\begin{itemize}
\item[(i)] If $a\neq  k\pi$ for any $k\in\mathbb N$, then
$$
u^* (x) \equiv \left[ {1\over a}\int_0^a  u_0^{1-p} \,{\rm d}x  \right]^{-{1\over p-1}}  \ \  \textrm{in}\  [0, a].
$$
\item[(ii)]  If $a =  k\pi$ for some $k\in\mathbb N$, then  one has the following statements.
    \begin{itemize}
    \item  When $1<p< 3/2$, denote
    \begin{equation}\label{thm-A0}
    A_0 = \left[ \frac{\int_0^{k\pi} (\cos x +1)^{1-p}\,{\rm d}x}{\int_0^{k\pi} u_0^{1-p} (x)\,{\rm d}x}\right]^{1\over p-1},
    \end{equation}
    then there exist $A$, $B_A$ with  $|A| \leq A_0$,  $B_A\geq |A|$ respectively such that
    $$
    u^*(x) = A\cos x +B_A  \ \  \textrm{in}\  [0,k\pi],
    $$
    where $B_A$  is the unique root of the equation
    \begin{equation}\label{thm-root-B}
    \int_0^{k\pi} (A\cos x +B_A)^{1-p} \,{\rm d}x = \int_0^{k\pi} u_0^{1-p} (x)\,{\rm d}x.
    \end{equation}
    In particular, if $|A| = A_0$,   then $u^* (x)= A \cos x + |A| $ in $[0, k\pi]$.
    \item  When $p \geq 3/2$,   there exist  $A$, $B_A$ with $|A|< B_A$ such that
    $$
    u^*(x) = A\cos x +B_A \ \  \textrm{in}\  [0,k\pi],
    $$
    where  $B_A$ is the unique root of the equation (\ref{thm-root-B}).
    \end{itemize}
\end{itemize}
\end{prop}



\begin{proof}
Recall that a typical property of this problem is as follows
\begin{equation}\label{conservation}
\int_0^a u^{1-p}\,{\rm d}x=\int_0^a u_0^{1-p}(x)\,{\rm d}x   \  \ \
  \forall \  \ 0\leq t<T_{{\max}}.
\end{equation}
Thus  according to Proposition \ref{prop-steadystates}(i), if $a\neq k\pi$ for any $k\in\mathbb N$, then
$$
u^*(x) \equiv \left[ {1\over a}\int_0^a  u_0^{1-p} \,{\rm d}x  \right]^{-{1\over p-1}} \ \ \textrm{in}\ [0,a].
$$
Proposition \ref{prop-steadystates-initial}(i) is proved.

It remains to discuss the case that  $a =  k\pi$ for some $k\in\mathbb N$. Thanks to Proposition \ref{prop-steadystates}(ii),  one sees that
$$
u^*(x) = A\cos x+ B \ \ \ \textrm{in}\ [0,k\pi]
$$
for some  $A,B\in  \mathbb R$   with  $  |A|\leq B$, where we include the possibility that $u^*$ is constant.
By (\ref{conservation}), one has
\begin{equation}\label{pf-conservation}
\int_0^{k\pi} (u^*)^{1-p}(x)  \,{\rm d}x = \int_0^{k\pi} (A\cos x+ B)^{1-p}  \,{\rm d}x =\int_0^{k\pi} u_0^{1-p}(x)\,{\rm d}x.
\end{equation}
It suffices to discuss the solvability of (\ref{pf-conservation}) with $ |A|\leq B$.
For clarity,  let us discuss the situations $1<p<3/2$ and $p\geq 3/2$ respectively.
\begin{itemize}
\item Assume that $1<p<3/2$,   then
$$
\int_0^{k\pi} (\cos x +1)^{1-p}\,{\rm d}x <+\infty.
$$
Thus  by (\ref{pf-conservation}), for $B\geq |A|$,
\begin{eqnarray*}
\int_0^{k\pi} u_0^{1-p}(x)\,{\rm d}x &=&\int_0^{k\pi} (A\cos x+ B)^{1-p}  \,{\rm d}x\\
&& \leq \int_0^{k\pi} (A\cos x+ |A|)^{1-p}  \,{\rm d}x= \int_0^{k\pi} (|A|\cos x+ |A|)^{1-p}  \,{\rm d}x,
\end{eqnarray*}
which yields that
$$
|A| \leq   \left[ \frac{\int_0^{k\pi} (\cos x +1)^{1-p}\,{\rm d}x}{\int_0^{k\pi} u_0^{1-p} (x)\,{\rm d}x}\right]^{1\over p-1} =  A_0.
$$
Thus  for give $0\leq |A|\leq A_0$,  on the one side, when $B=|A|$,
$$
\int_0^{k\pi} (A\cos x+ |A|)^{1-p}  \,{\rm d}x \geq \int_0^{k\pi} u_0^{1-p}(x)\,{\rm d}x.
$$
On the other side, it follows immediately from (\ref{pf-conservation}) that
\begin{equation}\label{pf-B-large}
 \lim_{B\rightarrow +\infty} \int_0^{k\pi} (A\cos x+ B)^{1-p}  \,{\rm d}x =0 < \int_0^{k\pi} u_0^{1-p}(x)\,{\rm d}x.
\end{equation}
Hence for give $0\leq |A|\leq A_0$, there exists a unique root in $[|A|, +\infty)$, denote by $B_A$, of the equation (\ref{thm-root-B}) as follows
\begin{equation*}
 \int_0^{k\pi} (A\cos x +B)^{1-p} \,{\rm d}x = \int_0^{k\pi} u_0^{1-p} (x)\,{\rm d}x.
\end{equation*}

\item  Assume that $p\geq 3/2$,   then
$$
\int_0^{k\pi} (\cos x +1)^{1-p}\,{\rm d}x = +\infty.
$$
Thus for any given $A\in \mathbb R$,
$$
\lim_{B\rightarrow |A|^+} \int_0^{k\pi} (A\cos x +B)^{1-p} \,{\rm d}x = \lim_{B\rightarrow |A|^+} \int_0^{k\pi} (|A|\cos x +B)^{1-p} \,{\rm d}x  =+\infty.
$$
This, together with (\ref{pf-B-large}), yields that there  exists a unique root in $(|A|, +\infty)$, denote by $B_A$, of the equation (\ref{thm-root-B}).
\end{itemize}

The proof  is complete.
\end{proof}

\section{Miscellaneous remarks}

Theorems \ref{thm-negativeenergy} and \ref{thm-nonnegativeenergy} are about the relations between the associated energy  functional $E[u](t)$ and the  finite-time  blowup/global existence  behavior of solutions to the problem (\ref{mainpb}). Theorem \ref{thm-negativeenergy} actually indicates that when $p>1$, if $E[u](t)$ becomes negative at some moment, then the solution will blow up in finite time. While Theorem \ref{thm-nonnegativeenergy}(ii) shows that when $1<p\leq 2$, if $E[u](t)$ always remains nonnegative, then the solution is global and bounded. An open question is what happens if $p>2$.

Theorem \ref{thm-convergence} itself verifies a fundamental property of the problem (\ref{mainpb}): global convergence of bounded solutions to some nonnegative steady state. In other words, $\omega[u_0] $ consists of a single steady state.

Furthermore, based on Theorem \ref{thm-convergence}, the classification of all nonnegative steady states obtained in Proposition \ref{prop-steadystates} and the integral conservation (\ref{integral}),
Proposition \ref{prop-steadystates-initial} reflects the influence of  the length of the interval $a$ and the parameter $p$ on the global dynamics of solutions to the problem (\ref{mainpb}). To be more specific, denote $\omega[u_0]  =\{ u^* \}$, we have the following picture.
\begin{enumerate}
\item If $a$ is not a multiple of $\pi$, then $u^*$ must be some positive constant uniquely determined by initial data.
\item If $a$ is a multiple of $\pi$ and $1<p<3/2$, due to the constrain of initial data, $u^*$  belongs to a continuum of nonnegative steady states and it is possible that $u^*$ touches zeros somewhere. In other words,  both {\it Case 1} and {\it Case 2} in the proof of Theorem \ref{thm-convergence} could happen.
\item If $a$ is a multiple of $\pi$ and $p\geq 3/2$, the constrain of initial data yields that $u^*$  belongs to a continuum of positive steady states, i.e., only {\it Case 1}  in the proof of Theorem \ref{thm-convergence} appears.
\end{enumerate}
It is also worth pointing out that in both (ii) and (iii), $u^*$  belongs to a continuum of steady states.  This partially explains why it is highly nontrivial to derive global convergence of solutions in Theorem \ref{thm-convergence}, especially when the problem (\ref{mainpb})   lacks comparison principles due to the introduction of nonlocal term.

\end{document}